\documentclass[conference]{IEEEtran}
\IEEEoverridecommandlockouts
\usepackage{cite}
\usepackage{amsmath,amssymb,amsfonts}
\usepackage{algorithmic}
\usepackage{graphicx}
\usepackage{textcomp}
\usepackage{xcolor}
\usepackage{booktabs}
\usepackage{multirow}
\usepackage[ruled,linesnumbered]{algorithm2e}
\usepackage{hyperref}
\usepackage[numbers,sort&compress]{natbib}

\hypersetup{hypertex=true,
colorlinks=true,
linkcolor=blue,
anchorcolor=blue,
citecolor=blue}

\def\BibTeX{{\rm B\kern-.05em{\sc i\kern-.025em b}\kern-.08em
    T\kern-.1667em\lower.7ex\hbox{E}\kern-.125emX}}

\begin{document}

\title{LRAMM —— Low precision approximates GEMM via RSVD
}

\author{
\IEEEauthorblockN{1\textsuperscript{st} 
Hongyaoxing Gu
}
\IEEEauthorblockA{\textit{
University of Chinese Academy of Sciences
} \\
BeiJing, China \\
guhongyaoxing23@mails.ucas.ac.cn
}

}

\maketitle

\begin{abstract}

Matrix multiplication computation acceleration has been a research hotspot across various domains. Due to the characteristics of some applications, approximate matrix multiplication can achieve significant performance improvements without losing much precision. 

In this paper, we propose LRAMM—a high-performance matrix multiplication approximation algorithm that combines mixed-precision quantized matrix multiplication with RSVD techniques, further enhancing efficiency within the error range of low-precision matrix multiplication by utilizing matrix low-rank decomposition technology. 
% Through evaluation, LRAMM can achieve a precision that is 10\% higher than int4 while simultaneously reducing the computational effort by 80\%(hand computing).
\end{abstract}


\section{Introduction}
Matrix multiplication is widely applied across various domains. In the field of deep learning, for instance, architectures such as GoogLeNet \cite{simonyan2014very}, ResNet \cite{he2016deep}, and VGG \cite{szegedy2015going} demonstrate that deeper networks often yield superior performance but require substantial training time to achieve high accuracy, exemplified by the 19-layer VGG network. It is well-known that approximately 90\% of the computational time in deep neural networks is spent on convolutional layers, which can be reformulated as matrix multiplications using techniques like img2col. In the realm of scientific computing, matrix multiplication also accounts for a significant proportion of time in applications such as molecular dynamics simulations and first-principles calculations. Therefore, optimizing the performance of matrix multiplication is crucial for enhancing the performance of these applications.

In response to this, numerous algorithmic improvements for matrix multiplication itself have been developed to accelerate the computation process. For example, Coppersmith and others\cite{coppersmith1987matrix} have developed an algorithm for matrix multiplication with a computational complexity of $O(n^2.375477)$, while Gall\cite{le2014powers} has developed an algorithm with a complexity of $O(n^2.3728639)$. Additionally, approximate algorithms for matrix multiplication leveraging the low-rank structure of matrices through methods such as singular value decomposition have been explored.

Furthermore, the low-rank structure of matrices has been proven to be of high value and ubiquity across numerous research fields. Applications leveraging the low-rank structure of matrices include but are not limited to imaging\cite{Lingala_Hu_DiBella_Jacob_2011}, fine-tuning large language models (\cite{Aghajanyan_Gupta_Zettlemoyer_2021},\cite{Wang_Tang_Duan_Wei_Huang_Ji_Cao_Jiang_Zhou_2021}), compressing neural networks\cite{Idelbayev_Carreira-Perpinan_2020}, and accelerating matrix multiplication\cite{Osawa_Sekiya_Naganuma_Yokota_2017}.

In the context of next-generation computing devices, high-dimensional matrices are represented by tensors, and low-precision operations carried out on tensors have garnered extensive support. Prominent among these heterogeneous platforms are Google's Tensor Processing Unit (TPU) \cite{jouppi2017datacenter}, Intel's Neural Network Processor (NNP), Neural Processing Units (NPU) \cite{hickmann2020intel,boutros2020beyond}, and NVIDIA GPU\cite{choquette2021nvidia}. These computing devices incorporate dedicated tensor units known as tensor cores\cite{Tensorcore}.

As indicated by official NVIDIA evaluations \cite{nvidia-ampere, turing-architecture, volta-architecture}, for FP16/FP32 mixed-precision deep learning, the A100 Tensor Core performance is \textbf{2.5} times that of the V100, increasing to \textbf{5} times with added sparsity. The peak speed for Int8 is \textbf{4} times that of FP16, and Int4 is \textbf{2} times that of Int8, while Int1 (binary) can be up to \textbf{4} times faster than Int4. Furthermore, the acceleration of Int8, Int4, and binary Tensor Core implementations supports deep learning inference, providing a robust foundation for the deployment of efficient and high-performance neural network models. These advancements in tensor core technology not only enhance computational throughput but also enable energy-efficient AI computing, paving the way for more scalable and accessible AI solutions across various industries and applications.

In this research, we introduce the use of low-rank approximation (RSVD) within the context of low-precision quantized matrix multiplication to enhance the efficiency of matrix multiplication operations. Our contributions are as follows:
\begin{itemize}
    \item Introduce several primary methods for matrix approximate multiplication and provide an analysis of their respective strengths and weaknesses.
    \item We propose an algorithm that integrates a mixed low-precision quantization strategy into low-rank approximate multiplication to improve computational speed.
    \item Conduct an analysis of computational error and efficiency for low-rank approximation and low-precision computation, and we design a series of experiments to demonstrate the efficacy of the proposed algorithm.
\end{itemize}

\section{Related works}

\subsection{Approximating matrix by SVD}
For a matrix $A$ of dimensions $m*k$, assume that it has rank of $p$. Then through Singular Value Decomposition (SVD), the matrix $A$ can be factorized into the product of three matrices:
$$ A = U \Sigma V^T  = \sum_{i=1}^p{\sigma_{i}u_iv_i^T} \in \mathbb{R}^{m*k}$$ 

where $ U $ is an $ m \times p $ orthogonal matrix, $ \Sigma $ is a $ p \times p $ diagonal matrix containing the singular values of $ A $, and $ V^T $ is the transpose of an $ k \times p $ orthogonal matrix $ V $. The diagonal entries of $ \Sigma $, known as singular values, are arranged in descending order, and each pair of singular values and corresponding columns from $ U $ and $ V $ represent an approximation of the original matrix $ A $ at varying levels of accuracy and rank reduction.$u_i,v_i$ is the row and column i in the matrix $U,V$.

Then, by retaining the first \( r \) non-zero singular values and their corresponding singular vectors in the \( U \) and \( V \) matrices, we can obtain an approximation \( A' \) of the original matrix \( A \). The \( r \) is chosen such that it captures the most significant information of \( A \) while reducing the complexity of the representation.

The approximate matrix \( A' \) can be computed as:

\[ A' = U_{r} \Sigma_{r} V_{r}^{T} = \sum_{i=1}^{r}{\sigma_{i}u_iv_i^T} \in \mathbb{R}^{m*k} \]

This approximation \( A' \) captures the essence of the original matrix \( A \) in a more computationally efficient form. By keeping only the top \( r \) singular values and their vectors, we are essentially performing a low-rank approximation and reducing the storage requirements while still maintaining a high level of accuracy, especially when the first few singular values account for a large portion of the matrix's energy or information content. This technique is widely used in numerical linear algebra and various applications, such as data compression, dimensionality reduction, and noise filtering, where a balance between accuracy and computational resources is crucial \cite{de2015data,anowar2021conceptual,paul2000transform}.

\subsection{Random SVD algorithms}
However, when $m$ and $m$ are both large, computing the full SVD of a matrix using numerical methods is typically computationally expensive and requires a significant amount of memory. As an effective alternative, randomized SVD approximation algorithms have recently garnered considerable interest and are competitive in computing fast low-rank approximations for large matrices\cite{zhou2014low,gupta2015deep}. The focus of randomized singular value decomposition algorithms is on effective sampling of significant matrix elements rather than neglecting the large matrix in a full SVD.

To approximate this computation using randomized algorithms,  here is a two-stage computation\cite{halko2011finding}:
\begin{enumerate}
    \item Compute an approximate basis for the range of $A \in \mathbb R^{m*n}$. We want a matrix $Q$ with $r$ orthonormal columns that captures the action of $A$. Formally, $A\approx QQ^{*}A$.
    \item Given such a matrix $Q$—which is much smaller than $A$— use it to compute our desired matrix factorization.
\end{enumerate}
In the case of RSVD, imagine we had access to $Q$. Then randomized SVD is the following:
\begin{algorithm}
    \label{algorithmofigemm}
  \SetAlgoLined
\KwData{$A,Q$(Intput matrix)}
  
  \KwResult{$U,\Sigma,V^{*}$(Svd decomposition matrix)}

    \tcc{Computes approximate matrix $B$}
    $\{B\}\leftarrow Q^{*}A$\;
    \tcc{Computes the svd of $B$}
    $\{U'\Sigma V^{*}\}\leftarrow B$\;
    $\{U\} \leftarrow QU'$\;

    \Return $U,\Sigma,V^{*}$\;
    
  \caption{A simple process of Randomized SVD }
\end{algorithm}

From the algorithm, it is evident that the computational requirements for RSVD are significantly reduced compared to the original SVD, with a time complexity of \( O(mnlog(r))+(m+n)r^2) \) - refer to \cite{halko2011finding}.

In order to find the orthogonal matrix $Q$ that RSVD can obtain higher precision. Numerous sampling strategies have been proposed, including uniform column/row sampling (with or without replacement)\cite{coppersmith1987matrix,le2014powers}, diagonal sampling or column norm sampling \cite{chellapilla2006high}, k-means clustering sampling\cite{chetlur2014cudnn}, and Gaussian sampling\cite{lin2014microsoft}. Consequently, compared to full SVD, randomized SVD methods offer higher memory efficiency and can typically achieve good low-rank approximations in a significantly faster manner.

\begin{figure}
	\centering 
\includegraphics[width=1\linewidth]{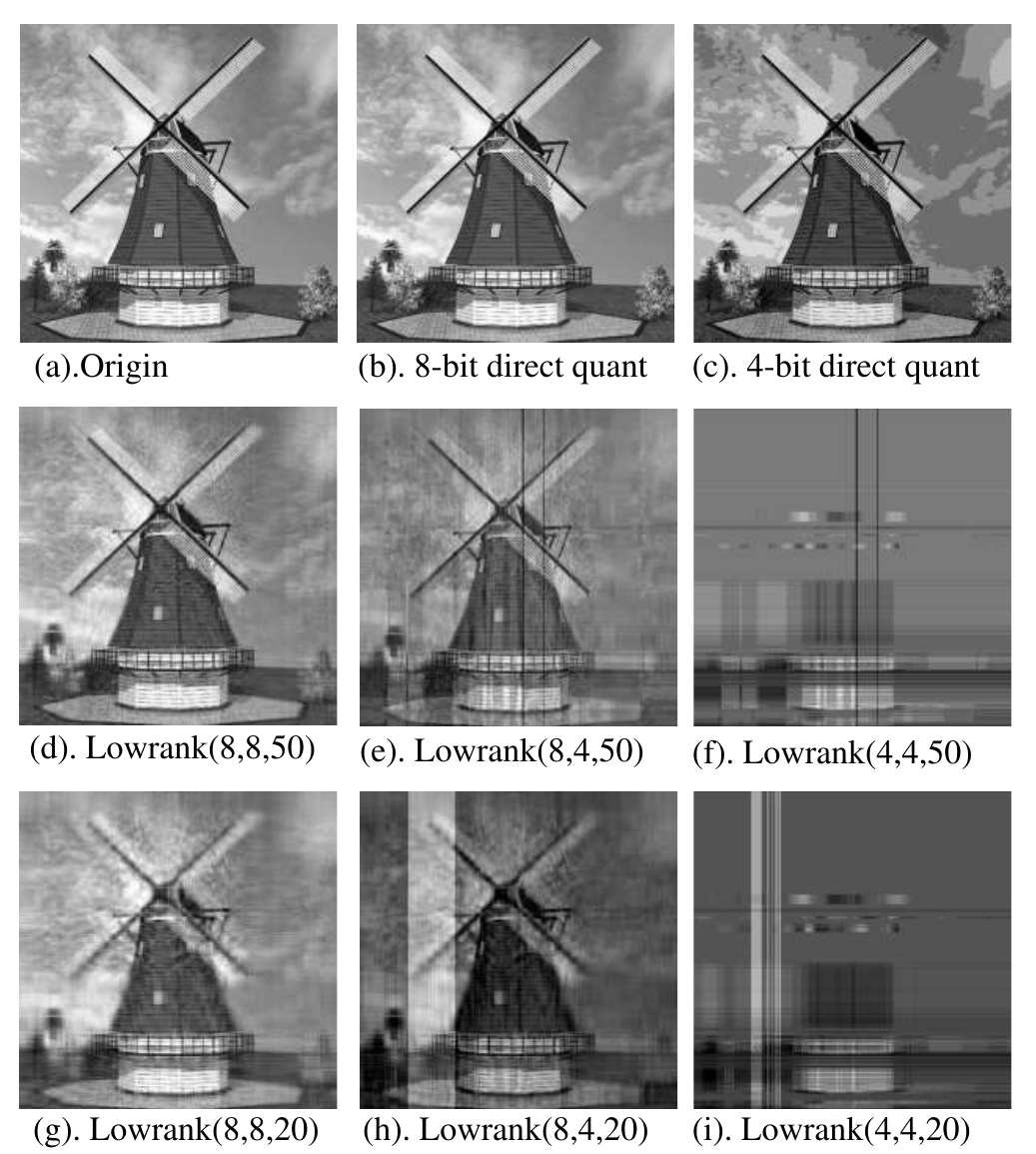}
	\caption{Using RSVD for low-rank approximation of images, where (a) represents the original image, and (b, c) are the full-size images using 8-bit and 4-bit quantization, respectively. (d-i) are images that use different approximation ranks and quantization bit-widths, denoted as Lowrank(d1, d2, r), where \( d_1 \) and \( d_2 \) represent the quantization bit-widths for the two matrices \( U\Sigma \) and \( V \), and \( r \) is the rank for the low-rank approximation.}
	\label{willow}
\end{figure}

Fig.\ref{willow} shows the effects of applying RSVD and low-precision quantization to an image.

\subsection{Approximate Matrix Multiplication}
Approximate Matrix Multiplication (AMM) is a technique designed to reduce the computational complexity of matrix multiplication operations. Matrix multiplication is a common operation across numerous scientific and engineering domains, particularly in the realms of machine learning and big data processing, where the scale and frequency of matrix computations are substantial. However, while traditional matrix multiplication algorithms, such as the Strassen algorithm, can theoretically reduce the amount of computation, they may encounter performance bottlenecks in practical applications, especially when dealing with large-scale datasets. The core idea of AMM is to approximate the results of matrix multiplication, thereby reducing the computational resources and time required. This approach typically involves several key aspects:

\textbf{Low-rank approximation:} This method involves decomposing a large matrix into the sum or product of several smaller matrices, which can be more easily multiplied. In this algorithm, the primary approach to dimensionality reduction of matrices is through the use of SVD (Singular Value Decomposition)\cite{Idelbayev_Carreira-Perpinan_2020}. However, the conventional SVD decomposition has a time complexity of $O(n^3)$, which is comparable to matrix multiplication and does not reduce execution speed. Therefore, the randomized SVD (RSVD) algorithm\cite{Drineas_Kannan_Mahoney_2006,yang2019robust,ji2016rank,drineas2006fast,eriksson2011importance} is usually applied in this method.

\textbf{Quantized Matrix Multiplication:} Quantization is the process of mapping continuous floating-point values to a finite set of discrete values. In the context of low-precision matrix multiplication, fixed-point numbers or quantized floating-point values are commonly utilized to approximate the elements of the original matrices. The choice of quantization strategy significantly impacts computational accuracy and efficiency. For instance, floating-point values can be quantized to integer values through techniques such as clustering\cite{equitz1989new} and scaling\cite{dai2021vs}. Additionally, methods involving KL divergence\cite{higham2019squeezing,adler2021quantization} entail truncating a portion of the original information before the mapping process, creating symmetric and well-distributed truncated information, which is then mapped into the integer domain.

\textbf{Matrix Multiplication Hashing Algorithm:} The core concept of this algorithm is to approximate the results of matrix multiplication using hashing functions and machine learning techniques, rather than performing standard multiply-accumulate operations. In the acceleration of neural networks, efforts have been made to utilize some form of hashing to speed up matrix multiplication\cite{spring2017scalable, chen2020slide}. 
Recently, machine learning is introduced into this process to learn the mapping relationships of matrix multiplication\cite{blalock2021multiplying}. Through extensive training data, the model can learn how to accurately predict the results of matrix multiplication. During the training process, the model receives the hash values of the input matrices and outputs the predicted multiplication results. By optimizing the parameters of the model, the algorithm can improve the accuracy of the predictions.

\section{Algorithm}
In the preceding section, we delved into the various methodologies for approximating matrix multiplication and the employment of low-precision computation through quantization in the context of matrix multiplication operations. In this section, we shall introduce an algorithmic approach for approximating matrix multiplication that utilizes both low-precision quantization and low-rank decomposition. Furthermore, in Section 3.2, we will present an analysis of the error associated with this algorithm, as well as an examination of its time complexity.

\subsection{Quantization GEMM}
\label{section:qgemm} 
We abbreviated quantization as $Q(A)$ and dequantization as $\widetilde{Q(A)}$, and employ a symmetric linear quantizer. Assume the representation with $N$ bits of integers is used, $a_{max}$ is the number with largest absolute value of matrix A.
\begin{equation}
\begin{split}
a_{int} = Q(a_{fp},\lambda) = TypeCast(\lambda*a_{fp},intN)
\label{Quantization}
\end{split}
\end{equation}

This function directly maps $a_{fp}$  proportionally to the range of integers. Where 
\begin{equation}
\begin{split}
\lambda=\frac{2^{N-1}-1}{a_{max}}
\label{get_lambda}
\end{split}
\end{equation}

As for the process of dequantization, which is the inverse of quantization, it involves restoring the results represented in integers back to floating-point numbers:
\begin{equation}
\begin{split}
a_{fp} = \widetilde{Q(a_{int},\lambda)} = TypeCast(a_{int}/\lambda,float)
\label{Dequantization}
\end{split}  
\end{equation}

For GEMM(general matrix multiplication), the quantized matrix multiplication process is shown in algorithm \ref{quantgemm}.
\begin{algorithm}
    \label{quantgemm}
  \SetAlgoLined
\KwData{$A,B,C$(Intput matrix);$\alpha$,$\beta$(Scalar);\newline N(Quant bit),$T_F$(Origin precision)}
  
  \KwResult{matrix D}

    \tcc{Computes quantized matrix}
    $\{A_{int},B_{int}\}\leftarrow Quant(\{A,B\},N)$\;
    \tcc{compute $D = A_{int}\cdot B_{int}$}
    $D_{int} =$ GEMM$(A_{int},B_{int})$\;
    \tcc{Computes dequantized matrix}
    $D_{F} \leftarrow \widetilde{Quant(D_{int},T_F)}$\;

    \Return $D = \alpha*D_{F} + \beta C$\;
    
  \caption{Algorithms of compute GEMM $D=\alpha A\cdot B+\beta C$ in Int-quantization }
\end{algorithm}

Assuming \( A \in \mathbb R^{m*k} \) and \( B \in \mathbb R^{k*n} \) with bit-widths $d_0$ , and quantization bit-widths is $d_1$. The Algorithms has a time complexity of $O\bigg(d_0^2(mk+kn+mn)+d_1^2(mnk)\bigg)$.

\subsection{Mixed low precision RSVD AMM}
\label{section:rsvd} 
Based on the RSVD approximation of the matrix form provided in the previous section, for the matrix multiplication \( C = A \* B \) where matrix $A$ of dimension $m*k$ and $B$ of $k*n$ , we can express the approximations of matrices \( A \) and \( B \) as:

\[ A \approx A' = U_{r} \Sigma_{r} V_{r}^{T}  \in \mathbb{R}^{m*k} \]
\[ B \approx B' = W_{r} \Gamma_{r} Z_{r}^{T}  \in \mathbb{R}^{k*n} \]

where \( A' \) and \( B' \) are the rank-\( r \) approximations of \( A \) and \( B \), respectively. \( U_{r} \) and \( W_{r} \) are orthogonal matrices, \( \Sigma_{r} \) and \( \Gamma_{r} \) are diagonal matrices containing the singular values and, in the case of \( B \), potentially other values depending on the specific decomposition, and \( V_{r} \) and \( Z_{r} \) are the matrices of right singular vectors.

To compute the product \( C' = A' \* B' \), we can leverage the properties of the SVD decomposition. The product of the approximate matrices \( A' \) and \( B' \) is given by:

\[ C' = (U_{r} \Sigma_{r} V_{r}^{T}) (W_{r} \Gamma_{r} Z_{r}^{T})  = (U_{r}\Sigma_{r})( V_{r}^{T} W_{r} )(\Gamma_{r} Z_{r}^{T}) \]

\begin{algorithm}
    \label{qrsvdgemm}
  \SetAlgoLined
\KwData{$A,B,C,D$(Intput matrix), $r$(ranf of low rank RSVD),$d_1,d_2,d_3$(bit-budget for each quantization GEMM)}
  
  \KwResult{$D$}

    \tcc{Computes RSVD approximate matrix of $A,B$}
    $\{U_{r},\Sigma_{r},V_{r}^{T}) \}\leftarrow RSVD(A,r)$\;
    $\{W_{r},\Gamma_{r},Z_{r}^{T}) \}\leftarrow RSVD(B,r)$\;
    \tcc{Calculate the product of the diagonal matrix $\Sigma_{r}, \Gamma_{r}.$}
    ${\widetilde{U_r},\widetilde{Z_r'}}\leftarrow GEMM( \{U_{r} \cdot \Sigma_{r},\Gamma_{r}\cdot Z_{r}^{T}\} )$\;
    \tcc{Calculate mix quantization GEMM}
    $E_{1}\leftarrow QM(V_{r}^{T}\cdot W_{r},d1)$\;
    $E_{2}\leftarrow QM(E_1\cdot \widetilde{Z_r'},d2)$\;
    $E_{3}\leftarrow QM(\widetilde{U_r}\cdot E_2,d3)$\;

    \Return  $D = \alpha*E_{3} + \beta C$\;
    
  \caption{Algorithms of mix quantization AMM implemented by RSVD.}
\end{algorithm}

The resulting approximation \( C' \) captures the most important information from both \( A \) and \( B \), and the multiplication is computationally more efficient due to the reduced rank. 

Subsequently, to further accelerate the computation outlined above, we will employ low-precision quantization operations. This technique involves reducing the bit-width of the numerical values involved in the calculations, which can significantly speed up the processing time while maintaining a reasonable level of accuracy.

Denote $U_{r}\Sigma_{r}$ as $\widetilde{U_r}$, $\Gamma_{r} Z_{r}^{T}$ as $\widetilde{Z_r}$. Then the approximate matrix multiplication can be expressed as $C = \widetilde{U_r}*V_{r}^{T}* W_{r}*\widetilde{Z_r}$.

Note that since matrices $\widetilde{U_r}$ and $\widetilde{Z_r}$ are obtained by multiplying diagonal matrices with orthogonal matrices, effectively, this operation is equivalent to scaling the orthogonal matrices $U_r$ and $Z_r$ by certain multiples along their rows and columns. The time complexity of this operation is $O(m*r) $, $ O(r*n)$, hence it can be directly executed using floating-point operations.

Consequently, three matrix multiplication operations remain. The quantization matrix multiplication is denoted as $C = QM(A*B,d)$, where d is the number of low precision integer quantization bits. Then the Low-precision quantized approximation of matrix multiplication can be formulated as follows:
$$C = QM(QM(QM(V_{r}^{T}* W_{r},d_1)*D1,d_2)*D2,d_3).$$

The pseudo code is provided in Algorithm \ref{qrsvdgemm}. Here, \( d_1 \), \( d_2 \), and \( d_3 \) represent the number of bits used for low-precision quantization in each step of the matrix multiplication process. The selection of different quantization bit widths can significantly impact the execution efficiency and accuracy of the algorithm, which will be discussed in a subsequent section.

\subsection{Approximation Error Analysis}
\label{section:error_analysis} 
For the convenience of error analysis, the following assumptions are made. Matrix $A \in \mathbb R^{m*k}$ with rank \( p \), and Matrix $B \in \mathbb R^{k*n}$ with rank \( q \). The rank used for the  RSVD approximated matrix is \( r \). The bit widths employed for the quantized multiplication are \( d_1 \), \( d_2 \), and \( d_3 \). 

As this algorithm employs a combination of randomized SVD decomposition, low-rank matrix approximation multiplication, and low-precision quantized matrix multiplication, analyzing the final error directly will be quite challenging due to the diverse approximation techniques involved. Our error analysis will separately examine the errors from these three aspects, and ultimately, derive the final error inequality of the algorithm by combining the errors from these components.

In this section, we will utilize the Frobenius norm $\Vert \cdot \Vert_F$ of matrices to conduct the analysis, where \(\mathbb E \) is represented as the mean.
\\
\\
\textbf{Errors in quantized matrix multiplication.}
For quantized matrix multiplication, $\widetilde{C}=\widetilde{A}*\widetilde{B}$ with $d_1,d_2$,  $R_1:(max(A))$, $R_2:(max(B))$, $\lambda_1 = \frac{2^{d_1-1}-1}{R_1}$,$\lambda_2 = \frac{2^{d_2-1}-1}{R_2}$(From equation.(\ref{get_lambda})). We have(Proof in Theorem.\ref{svdfinq})
\begin{equation}
\label{er_quant}
\begin{split}
\mathbb E(\Vert C-\widetilde{C} \Vert_F) &\leq k(\sigma_1\lambda_2^{-1}\sqrt{n}+\gamma_1\lambda_1^{-1}\sqrt{m} \\
    &+\lambda_1^{-1}\lambda_2^{-1}\sqrt{mn})
\end{split}
\end{equation}
Besides 
$$\mathbb E(\Vert C\Vert_F) = \sqrt{\sum_m\sum_n(\sum_ka_{ik}*b_{kj})^2}\leq\sqrt{mnk}R_1R_2$$
Use the ratio of the Frobenius norm of the error matrix to that of the original matrix as the measure of relative error, denoted as \( E_r \). The relative error for this process can be expressed as:
\begin{equation}
\label{er_qgemm}
\begin{split}
\mathbb E(E_r) &= \frac{\mathbb E(\Vert C-\widetilde{C} \Vert_F)}{\mathbb E(\Vert C\Vert_F)} 
    = \sqrt{\frac{k}{m}}\frac{\sigma_1\lambda_2^{-1}}{R_1R_2} \\
    &+ \sqrt{\frac{k}{n}}\frac{\gamma_1\lambda_1^{-1}}{R_1R_2} + \frac{\lambda_1^{-1}\lambda_2^{-1}}{\sqrt{k}R_1R_2}
\end{split}
\end{equation}
Notice that the last term is much smaller than the first two and can be ignored. And m,n,k, as the dimensions of the matrix, can be regarded as the same order of magnitude. It can be deduced that
\begin{equation}
\label{er_qgemm}
\begin{split}
\mathbb E(E_r) \propto \{2^{d_1},2^{d_2}\}
\end{split}
\end{equation}
That is to say, once the quantization bit widths \( d_1 \) and \( d_2 \) are established, the quantized matrix multiplication is invariant to \( m \), \( n \), and \( k \), and is solely dependent on the distribution of the matrix values. This finding aligns with results from actual testing in Fig.\ref{merge_precision}.
\\
\\
\textbf{Errors in low-rank matrix approximation via RSVD.}

For the matrix $A$ has SVD decomposition. If we approximate matrix \( A \) using a submatrix of rank \( r \), denote $p= min(m,n)$
$$A_r = \sum_{i=1}^r{\sigma_{i}u_iv_i^T} \in \mathbb{R}^{m*k} \  , \ (r \leq p)$$ 
then the error at this stage can be upper bounded as (Proof in Theorem.\ref{svdfinq}):
\begin{equation}
\label{svdfinq_z}
\begin{split}
\mathbb E||A - A_r||_F \leq \sigma_{r+1}\sqrt{p-r} 
\end{split}
\end{equation}

However, this represents an ideal scenario. Complete SVD decomposition requires the computation of all eigenvalues of the matrix, which entails a significant amount of redundant computation compared to the submatrix we require, with a time complexity of \( O(n^3) \). Therefore, for our AMM, we resort to the RSVD decomposition which will save a considerable amount of computational effort and introduce some errors. The error associated with the approximated matrix using this method is as follows(Proof in Theorem.\ref{rsvderror}):
\begin{equation}
\label{rsvderror_z}
\begin{split}
\mathbb E \Vert A-A_r\Vert \leq \left[ 1+4\sqrt{\frac{2p}{r-1} }\right ]^{1/(2q+1)}\sigma_{r+1}.
\end{split}
\end{equation}

Here, \( q \) is a parameter included in RSVD, representing the number of iterative solutions. A higher \( q \) leads to more accurate computational results. In this algorithm, \( q=0 \) is chosen. Given that the algorithm involves low-precision computations, high precision for RSVD is not required.

Due to the characteristics of the RSVD algorithm, this form of error is more precise than the direct truncation of SVD. Since \( r < p \), the error of RSVD will not exceed $7\sigma_{r+1}$. From these two forms, it can be observed that the RSVD algorithm.

Compared to the errors introduced by low-precision quantization, it can be stated that, in most cases, the RSVD does not introduce a higher error term compared to the original SVD.
\\
\\
\textbf{Errors in low-rank RSVD AMM(LRAMM).}
In the previous section, we analyzed the general error case of quantized AMM. However, the form presented contained too many parameters, which was not conducive to intuitive analysis. In this section, we will impose certain restrictions on matrices \( A \) and \( B \), as well as the quantization conditions, to provide a more intuitive form of error.

Assuming the matrix dimensions are \( m = n = k \), and that matrices \( A \) and \( B \) are identically distributed with a maximum value of \( R \). 

Denote $2^{d_i-1}-1$ as $D_i$,  $\frac{\sigma_{r+1}^2}{\sigma_1^2}(k-r)$ as $f(r)$ which means under a fixed distribution for matrices, the magnitude of this term depends solely on the value chosen for \( r \).
Then(Proof in Theorem.\ref{spec_amm_1})
\begin{equation}
\label{spec_amm_0}
\begin{split}
\mathbb E(\Vert {C'_r} - C \Vert_F) &\leq k\sigma_1^2\sqrt{r} \bigg(\sqrt{-f(r)+\frac{1}{D_1}+\frac{1}{D_2}+\frac{K}{D_1D_2}}  \\
        &+ \sqrt{f(r)+\frac{1}{D_2}+\frac{1}{D_3}+\frac{K}{D_2D_3}} \bigg ) \\
        &+kr\sigma_1^2\bigg( \sqrt{f(r)+\frac{1}{D_1}+\frac{1}{D_2}+\frac{K}{D_1D_2}} \bigg)\\
        &\cdot \bigg ( \sqrt{f(r)+\frac{1}{D_2}+\frac{1}{D_3}+\frac{K}{D_2D_3}} \bigg)
\end{split}
\end{equation}
Although this formulation appears complex at first glance, analysis can reveal three key takeaways:
\begin{itemize}
    \item For a given matrix, its largest singular value \( \sigma_1, \gamma_1 \), and dimensions \( m \), \( n \), \( k \) remain constant. These factors are irrelevant to the final relative error \( E_R \).
    \item The target rank \( r \) for the final relative error is not linear but is related to both the rank of the original matrix and the number of quantization bits. It is not the case that selecting a higher rank \( r \) for the approximate matrix will always lead to higher precision in the algorithm's final output(This point may seem counterintuitive, but subsequent experimental sections will substantiate the correctness of this observation.).
    \item The impact of the quantization bit widths \( d_1 \), \( d_2 \), and \( d_3 \) on the algorithm's error is not same. From the error formula(\ref{spec_amm_0}), \( d_2 \) has a greater influence on the error than the other two bit widths.
\end{itemize}

\subsection{Time complexity analysis}
\label{section:complexityr_analysis} 
To account for the efficiency of low-precision computation, we incorporate the input bit-width into the time complexity analysis. Assuming the original input precision is \( d_0 \) bits, and the low-precision quantization bit-widths are \( d_1 \), \( d_2 \), and \( d_3 \), respectively.

For LRAMM (Low-Rank Approximate Matrix Multiplication), the algorithm primarily consists of three components. As shown in the algorithm.\ref{qrsvdgemm}:
\begin{enumerate}
    \item Perform an RSVD decomposition with rank \( r \) on the two original matrices \( A \) and \( B \), which is of complexity \( O(mklog(r))+(m+k)r^2)) \) and  \( O(nklog(r))+(n+k)r^2)) \). A total of 
    \begin{equation}
    \label{time_c_1}
    \begin{split}
     O\bigg(d_0^2\big((m+n)klog(r)+(m+n+2k)r^2)\big)\bigg) 
    \end{split}
    \end{equation}
    
    \item Compute \( \widetilde{U_r} \) and \( \widetilde{Z_r} \) by performing diagonal matrix multiplication at the original precision, the time complexity of this step is 
    \begin{equation}
    \label{time_c_2}
    \begin{split}
     O\bigg( d_0^2(mr+nr) \bigg)
    \end{split}
    \end{equation}

    \item Computes three low-precision matrix multiplications ($E_1,E_2,E_3$--algo.\ref{qrsvdgemm}). The time complexity is respectively $ O\bigg(2d_o^2kr+d_2^2kr^2\bigg) $, $O\bigg(d_o^2(r^2+nr))+d_3^2nr^2\bigg)$, $O\bigg(d_o^2(nr+mr))+d_1^2mnr\bigg)$. A total of 
    \begin{equation}
    \label{time_c_3}
    \begin{split}
     O\bigg (d_o^2r(2k+r+2n+m)+ r(d_1^2mn+ d_2^2kr+d_3^2nr)\bigg) 
    \end{split}
    \end{equation}

\end{enumerate}

To intuitively demonstrate the effectiveness of the algorithm, we assume \( m = n = k \) and \( d_1 = d_2 = d_3 \). Under these conditions, the complexity of the algorithm is:
\begin{equation}
\label{time_c_4}
\begin{split}
 O\bigg( d_0^2(2m^2log(r)+ 4mr^2+7mr+r^2)+d_1^2(m^2r+2mr^2)\bigg) 
\end{split}
\end{equation}    

From this expression, it is shown that when \( r \) is small, the computational effort for $d_0^2$ part will not be significant. The primary computational load lies within the low-precision matrix multiplication, and in practical computations, the computers support only a limited number of bits for low-precision integer widths, typically 1, 4, and 8 bits, while origin precision is at least 32.

Compared to the original matrix multiplication(about $O(d_0m^3)$), our algorithm exhibits better efficiency. When contrasted with low-precision matrix multiplication(about $O(d_1m^3)$), our algorithm also manages to save a considerable order of magnitude in computational effort while $r$ is small.

\subsection{Parameter turning}
 In LRAMM there are four parameters that need to be specified: the rank \( r \) of the approximate matrix; and the quantization bit widths \( d_1 \), \( d_2 \), \( d_3 \) for the three matrix multiplications. Below, we will discuss the relationship between the selection of these parameters and how it determines the precision and efficiency of the algorithm.

In \ref{section:complexityr_analysis}, through the final provided expression of time complexity, it is evident that the algorithm's runtime \( Time \propto \{r^{-1},d_i^{-1}\} \), indicating that opting for lower quantization bit widths and a lower approximate rank will result in a reduced computation time.

From the error analysis section in \ref{section:error_analysis}, the error has an expression that is analogous to a quadratic form with respect to \( r \). In the actual algorithm, the relationship between the size of \( r \) and computational precision is akin to the graph of a quadratic function, and it also yields good results when it approaches 1/10 of the matrix dimensions. This will be demonstrated in the subsequent experimental section \textcolor{red}{need ref}.

Regarding the selection of \( d_1 \), \( d_2 \), and \( d_3 \), since 1-bit quantization can lead to a significant amount of computational error, and most systems only support the selection of 4 or 8-bit integers, it is advisable to choose \( d_2=4 \) as and $d_1,d_3=8$ to achieve better general effect, as analyzed in the previous section. More detailed results can be found in the experimental section \textcolor{red}{need ref}.

\section{Evaluation}

To evaluate the effectiveness of LRAMM, we conducted experimental tests from multiple aspects, including precision testing under different scales and distributions, as well as efficiency testing of the algorithm's execution time.

The experimental setup involved the following computing devices:
\begin{itemize}
    \item Model: Intel(R) Xeon(R) Platinum 8163 CPU
    \item Compilation:  \texttt{-O3 -mavx,-mfma,-march=native}
    \item Operating System: Ubuntu 11.4.0-1ubuntu1~22.04
\end{itemize}

The algorithm implementation, as described in Section 3 , consists of three main components:

\begin{itemize}

    \item \textbf{RSVD}: \texttt{Eigen} on CPU was utilized to construct the RSVD kernel.
    \item \textbf{GEMM}: \texttt{Eigen} on CPU was employed for matrix multiplication. 
    \item \textbf{Matrix Quantization and Reduction Operations}: Developed using the \texttt{OpenMP} on CPU for matrix quantization and reduce operations.
    
\end{itemize}

\begin{figure}
	\centering
\includegraphics[width=1\linewidth]{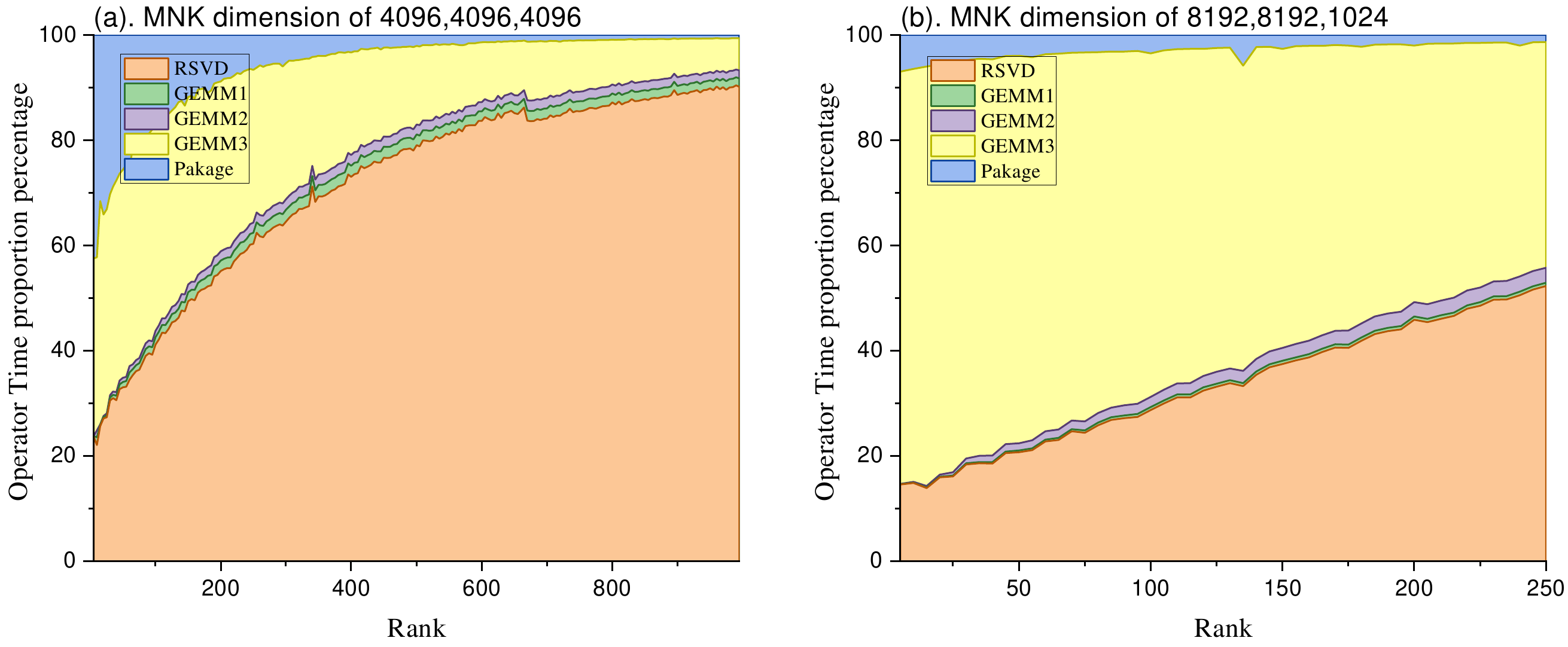}
	\caption{The proportion of operator execution time in LARMM across different scales, where GEMM1$ \sim $ 3 denote the running times of the three low-precision matrix multiplications in Algorithm 3, RSVD represents the running time for randomized SVD, and package refers to the running time for quantization operations and other matrix operations.}
	\label{time proportion test}
\end{figure}

\subsection{Time proportion test}
Algorithm execution time proportion test. The execution time proportion of the main parts of the algorithm is shown in Fig.\ref{time proportion test}, with a matrix scale of \( 4000^3 \), where "other" includes operations such as quantization and matrix addition.

It can be observed that when the rank is small, the execution times of the three are comparable. However, when the approximate rank is relatively large, the main execution time of the algorithm is concentrated on RSVD and matrix multiplication. In practical use, the approximate SVD matrix algorithm will only be applied when the matrix rank is not large (one fifth of the matrix dimension or less). 

Additionally, among the three matrix multiplication items, it can be seen that the time taken by GEMM3 is much greater than that of GEMM1 and GEMM2. From the error analysis presented earlier, it is evident that using a higher precision for GEMM1 compared to the other two multiplication steps can achieve better accuracy, which also indicates that the performance improvement brought by using a lower precision (int4) for the third matrix multiplication is quite significant, confirming the effectiveness of the variable-bit-width quantization algorithm.

\begin{figure}
	\centering
\includegraphics[width=1\linewidth]{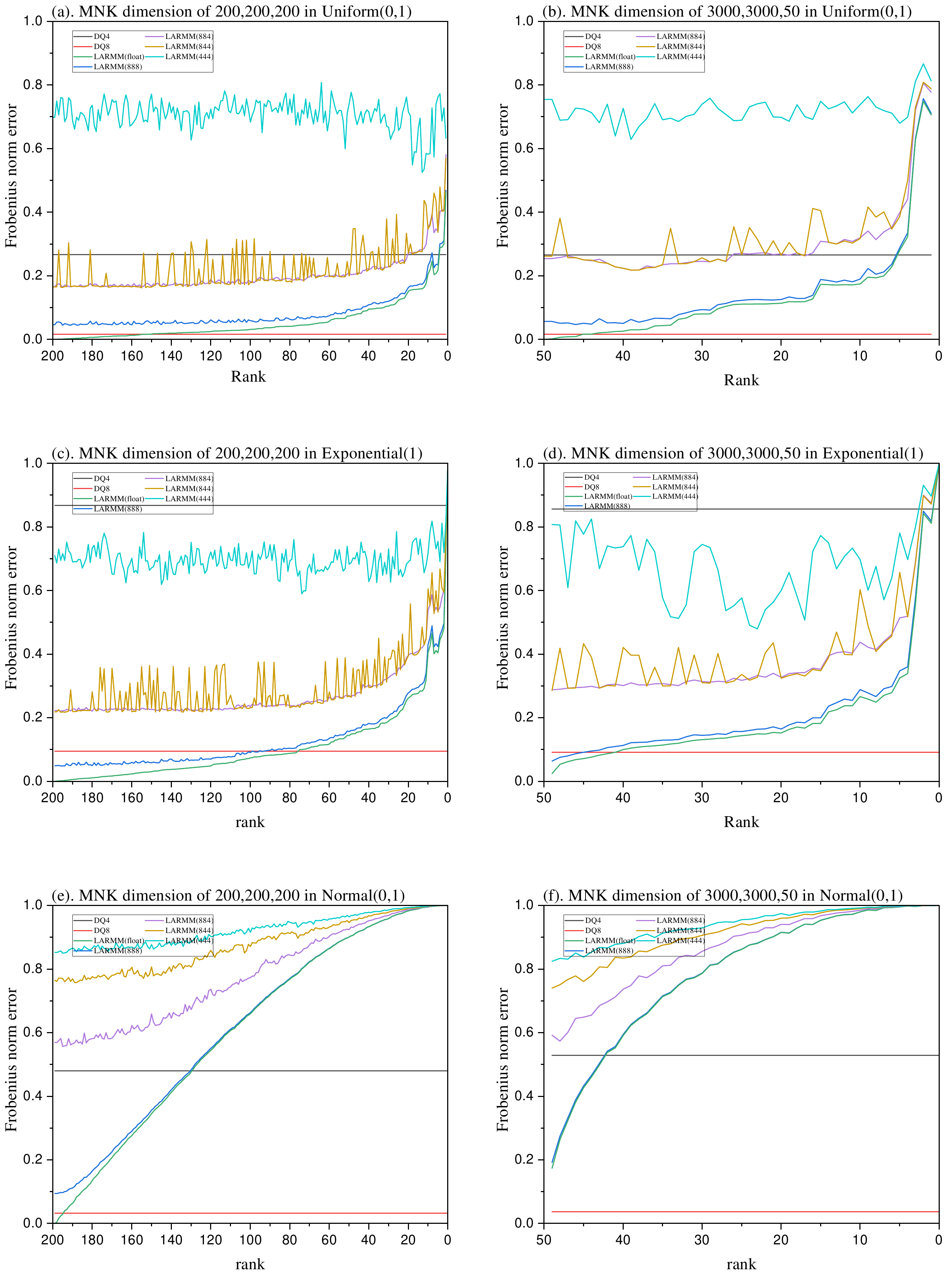}

	\caption{The relative error and the approximate rank of LARMM with different parameters under matrices of various distributions. Where \( DQ4 \) and \( DQ8 \) denote the fully-sized quantized matrix multiplication with 4-bit and 8-bit precision, respectively. \( LARMM(d_1d_2d_3) \) indicates that the LARMM algorithm uses low-precision quantization bit-widths \( d_1 \), \( d_2 \), \( d_3 \) for the three steps of matrix multiplication. \( LARMM(float) \) signifies that the matrix multiplication in LARMM is directly computed using single precision without employing low-precision quantization. }
	\label{merge_precision}
\end{figure}

\subsection{Precision test}
In the precision tests, we used matrices of different scales and distributions for testing as shown in Fig.\ref{merge_precision}. It can be observed from the figure that the algorithm's error is consistent with the analysis in Section 3. The error is subject to multiple influencing factors: for a specific distribution, the relative error of direct quantization is often a fixed number, which is smaller when the matrix is more uniformly distributed and tends to be larger under non-uniform distributions, such as when there are extreme values (exponential distribution).

Regarding the selection of quantization bit widths \( d_1 \), \( d_2 \), \( d_3 \), the results of the algorithm align with the analysis from the previous section. The relative error for using higher computational precision for \( d_1 \) and \( d_2 \) (LARMM884, LARMM844) is significantly smaller than that of LARMM444. This implies that employing higher precision for the two matrix multiplications with the shortest computation time can lead to a notable improvement in accuracy. Moreover, it can be observed that under both exponential and uniform distributions, LARMM884 and LARMM844 can achieve a precision better than int4 by using an approximation rank that is about 1/10th of the original matrix rank. At the same time, since the computational effort is concentrated on matrix multiplication, the algorithm can also achieve a significant speedup ratio. A detailed analysis can be found in the following section.

However, for the normal distribution, the performance of the matrix approximation algorithm is not as effective. As shown in Fig.\ref{merge_precision}(e,f), the algorithm requires approximately half of the rank to achieve the same accuracy as int4. The results indicate that direct quantization is a better approach. This is due to the nature of low-rank matrix approximation multiplication. It is important to recognize that compression techniques based on low-rank decomposition are only effective when the matrix being compressed is inherently low-rank, implying that \(\|A_r - A\|_F\) is initially small. When \(\|A_r - A\|_F \approx 0\), LARMM can achieve the same error level as direct quantization with less computational effort.

\begin{figure}
	\centering
\includegraphics[width=1\linewidth]{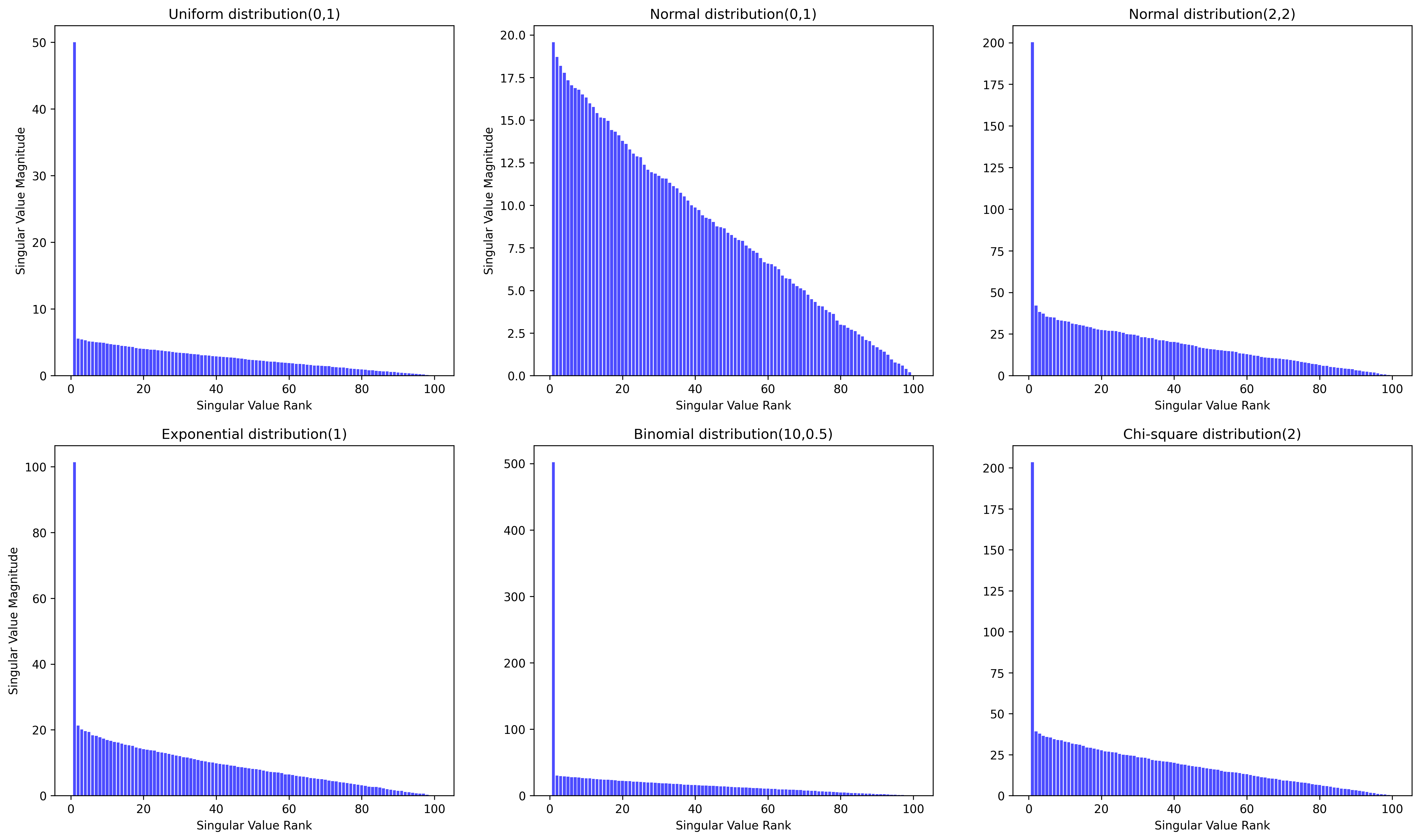}
	\caption{The distribution of singular values in matrices under different distributions, where the matrix dimensions are \( 100 \times 100 \).}
	\label{singular_value_bar_charts}
\end{figure}

Fig.\ref{singular_value_bar_charts} shows the distribution of singular values for matrices under different distributions. When a matrix is non-negative, its rank distribution is often more uneven, with some singular values that are significantly larger than the others. This means that for such matrices $\sigma_r/\sigma_1 \approx 0$ then \(\|A_r - A\|_F \approx 0\), and low-rank approximation algorithms can often achieve better results. In contrast, for matrices with more uniformly distributed singular values, such as those following a Normal(0, 1) distribution, this condition is not met, and thus the use of low-rank approximation algorithms is often not suitable.

% Here, two sets of test results are provided for different matrix scales, where the matrices used are uniformly distributed matrices. And from which it can be observed that 
% \begin{enumerate}
%     \item the relative performance of the 884 configuration is quite good.
%     \item The overall error of the algorithm is quite impressive. At a rank that is 1/10th of the original, the precision of algorithm 884 surpasses that of int4.
%     \item Due to the use of a low-rank approximation with \( r \), and compared to full int4 quantization, although algorithm 884 includes two additional low-precision matrix multiplications with a precision of int8, the complexity of each of these multiplications is \( O(kr^2) \). Since \( r \) is quite small, the actual computational effort is minimal compared to the last int4 multiplication. As a result, the overall performance is superior to that of full int4 quantization.
% \end{enumerate}

\begin{figure}
	\centering
\includegraphics[width=1\linewidth]{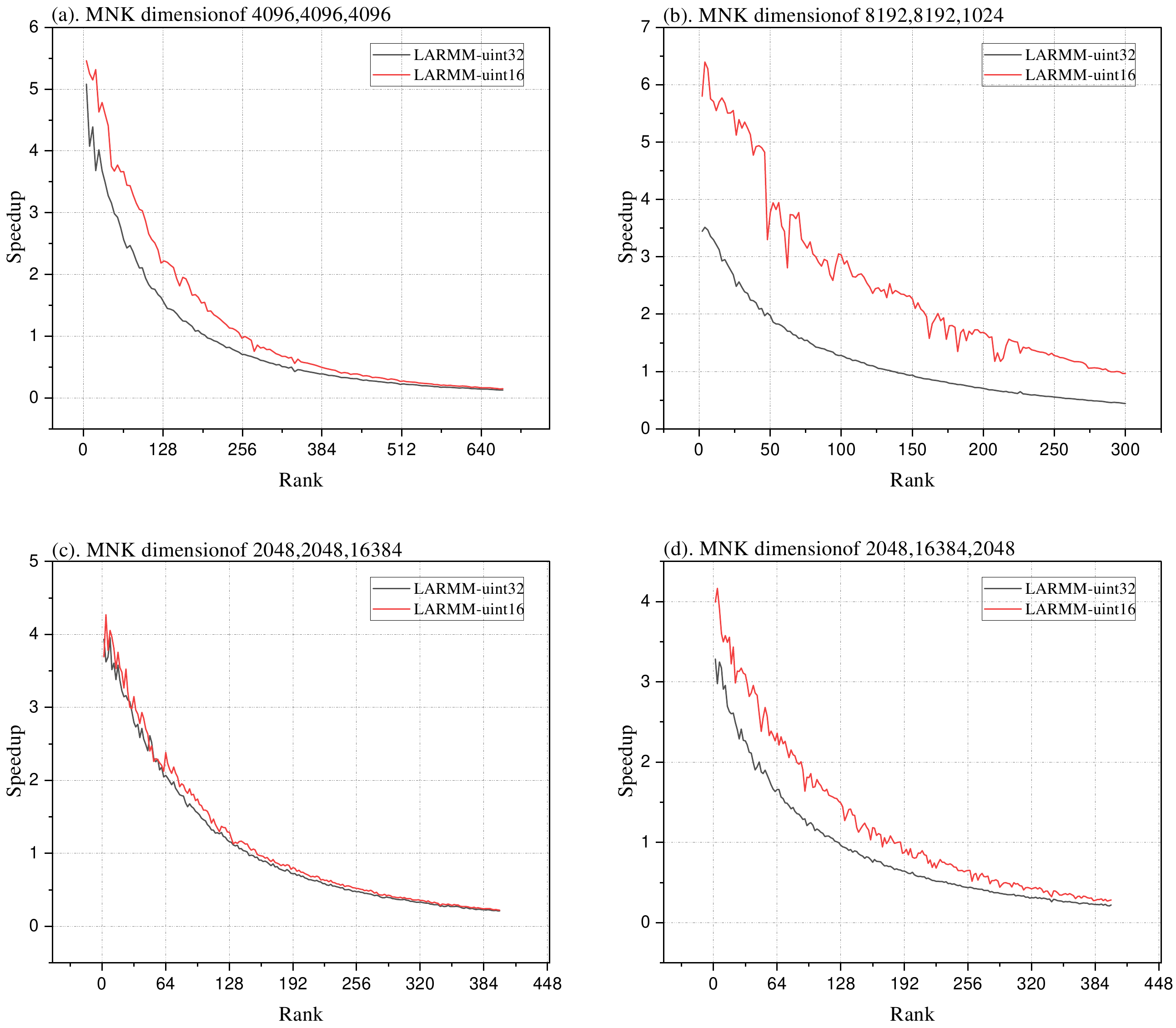}
	\caption{The acceleration ratio of LARMM with different bits in different approximate rank. The full quantization matrix multiplication using uint16 as the baseline, where LARMM-uint32 indicates that the matrix multiplication employs 32-bit computation, and LARMM-uint16 signifies that the matrix multiplication uses 16-bit computation.
}
	\label{merge_speed}
\end{figure}
\subsection{Algorithm efficiency test}
Due to the fact that MKL (Math Kernel Library) only supports 16-bit and 8-bit low-precision integer matrix multiplications, and the 8-bit matrix multiplication has input sign limitations that are not conducive to use, in this part of the testing, the full quantization matrix multiplication using uint16 is used as the baseline, representing direct quantization computation at low precision. Fig.\ref{merge_speed} presents the speedup ratio of the LARMM algorithm under different matrix scales and different operand bit-widths. 

Compared to full-size low-precision quantized matrix multiplication, LARMM can achieve significant speedup ratios within a certain approximate rank, and performs noticeably better when the matrix is of low rank. Moreover, when \( K \) is relatively small, the low-precision LARMM performs better. In the case of matrices with dimensions 8192×8192×1024, the performance of the 16-bit LARMM is nearly twice as good as that of the 32-bit version, which also demonstrates the effectiveness of using low-precision computation for intermediate GEMM (General Matrix Multiply) operations.

Combining the precision analysis from the previous section, it is often possible to represent matrices of certain distributions using a low-rank matrix corresponding to only the first few singular values, and this size is independent of the matrix scale. This implies that LARMM can achieve excellent results with matrices of some distributions. For instance, in a \( 8192 \times 8192 \times 1024 \) binary uniform matrix, using an approximate matrix multiplication with a rank of \textbf{50} can yield an acceleration ratio nearly \textbf{4} times that of the original low-precision quantization method.

\section{Discussion}
\textbf{Fusion operation: }Currently, LARMM uses separate RSVD and matrix multiplication quantization operations, which require re-quantization and computation based on the bit-width of the operation after each calculation. These steps can be fused to enhance computational efficiency. 

\textbf{Low-precision RSVD algorithm:} Since the algorithm adopts low-precision for key steps and the RSVD step occupies a significant portion of the algorithm's execution, using an even lower precision for computation during the RSVD step may be a good option. A mixed-precision matrix low-rank computation method has been proposed in LPLR\cite{saha2023matrix}.

\textbf{Computation on GPUs:} Low-precision operations are better supported on GPUs, and compared to CPUs, low-precision computation on GPUs can achieve a more significant speedup ratio. However, due to the limited availability of RSVD algorithms on GPUs\cite{RSVDPACK}, and the fact that existing algorithms were developed earlier and do not support low-precision computation, it would be beneficial to develop efficient mixed-precision RSVD algorithms tailored to the latest GPU architectures.

\textbf{Auto rank:}

\section{Conclusion}
In this work, we have presented a low-precision low-rank approximate matrix multiplication algorithm that employs mixed low-precision computation in key matrix multiplication operations to enhance computational efficiency. An error analysis of the algorithm is conducted, and based on this analysis, a rational strategy for the selection of quantization parameters is provided. The validity of the theoretical analysis and the effectiveness of the algorithm are then demonstrated through experiments.

\clearpage

% \begin{table*}
% \caption{\textbf{Mathematical symbols used in paper}}
% \centering
% \begin{tabular}{ccc}%四个c代表该表一共四列，内容全部居中
% \toprule
% Notation & Description & Comment\\
% \midrule%第一道横线
% \multirow{4}{*}{A,B}& \multirow{4}{*}{origin matrix}& {If matrix \( A \) appears alone, its dimensions}\\%表格宽度参数采用*代表自动宽度
% && are \( m \times n \). If both \( A \) and \( B \) are present\\
% && , then the dimensions of \( A \) are \( m \times k \) \\
% && and \( B \) are \( k \times n \).\\
% \midrule%第二道横线 
% \multirow{2}{*}{m,n,k}&Dimensions   & $A\in \mathbb R^{m*n} or R^{m*k}, $ \\
% &of input matrix & $B\in \mathbb R^{k*n} $\\

% \midrule
% \multirow{2}{*}{$A = U\Sigma V^T$}&Svd of $A$   & $ If A\in \mathbb R^{m*n}, U\in \mathbb R^{m*m}$ \\
% && $B\in \mathbb R^{k*n} $\\

% \midrule%第三道横线 
% Data5&Data6&Data7 \\
% \bottomrule%第四道横线
% \end{tabular}
% \end{table*}

\begin{table*}
\caption{\textbf{Mathematical symbols used in paper}}
\centering
\begin{tabular}{ccc}%四个c代表该表一共四列，内容全部居中
\toprule
Notation & Description & Comment\\

\midrule%第二道横线 
\multirow{1}{*}{m,n,k}&Dimensions of input matrix  & $A\in \mathbb R^{m*n} or R^{m*k}, B\in \mathbb R^{k*n} $ \\

\midrule
$p$ & the min one in matrix dimensions  & $p=min(m,n)$\\
\midrule
$r$ & the rank used in low rank approximation  & none \\

\midrule%第一道横线
\multirow{2}{*}{A,B}& \multirow{2}{*}{origin matrix}& {If matrix \( A \) appears alone, its dimensions are \( m \times n \). If both \( A \) and \( B \) }
\\ &&{are present , then the dimensions of \( A \) are \( m \times k \)  and \( B \) are \( k \times n \).}\\

\midrule
\multirow{1}{*}{$M$}&Max absolute number in matrix& $M = max(a_{ij})$ \\

\midrule
\multirow{1}{*}{$U\Sigma V^T$}&Svd of $A$   & $ If A\in \mathbb R^{m*n}, U\in \mathbb R^{m*m}, \Sigma \in \mathbb R^{p*p},  V\in \mathbb R^{n*n}$ \\

\midrule
\multirow{1}{*}{$W\Gamma Z^{T} $}&Svd of $B$   & $ If B\in \mathbb R^{m*n}, W\in \mathbb R^{m*m}, \Gamma \in \mathbb R^{p*p},  Z\in \mathbb R^{n*n}$ \\

\midrule
\multirow{1}{*}{$u, v, z, w$}&A column in the matrix $U,V,Z,W$& $A = \sum_{i=1}^p{\sigma_{i}u_iv_i^T} \in \mathbb{R}^{m*n} $
\\

\midrule
\multirow{1}{*}{$\sigma,\gamma$}&The eigenvalue of the matrix & In this paper $\sigma_i$ is the ith eigenvalue of $A$ and $\gamma_i$ of $B$  \\

\midrule
\multirow{1}{*}{$A_r$}&Low rank approximation with rank r   & $A_r = U_r\Sigma_r V^T_r$ \\

\midrule
\multirow{1}{*}{$\mathbb E()$}&Expected Value symbol &none  \\

\midrule
\multirow{1}{*}{$\Vert\cdot\Vert$}&Norm of matrix & $\Vert\cdot\Vert_2 = \sigma_1 \ , \ \Vert\cdot\Vert_F = \sqrt{\sigma_1 + \dots + \sigma_p}$  \\

\midrule
\multirow{1}{*}{$d$} & The bit-width used in quantized matrix multiplication & The three quantized gemm in this paper are $d_1, d_2, d_3$ respectively  \\

\midrule
\multirow{1}{*}{$\lambda$} & The parameters used for the quantization operation & $\lambda = \frac{2^{d-1}-1}{M}$  \\

\midrule
\multirow{2}{*}{$R_A,R_B$} & Represents the difference matrix between the & \multirow{2}{*}{eg.$R_A = A-\widetilde{A}$}  \\
&approximate matrix and the original matrix.&\\

\midrule
\multirow{1}{*}{$D$} & A variable used to simplify an expression & $D_i = 2^{d_i-1}-1$  \\

\bottomrule%第四道横线
\end{tabular}
\end{table*}

\clearpage

\textbf{Appendix}
% \begin{tabular}{cccccc}
%    \toprule
%     Notation & Description & Comment  \\
%    \midrule
%    \multirow{2}*{A,B} & {origin matrix} & When matrix \( A \) appears alone, its dimensions are \( m \times n \). If both \( A \) and \( B \) are present, then the dimensions of \( A \) are \( m \times k \) and \( B \) are \( k \times n \).  \\
%     ES & F & 15  \\
%     SE & M & 17 \\
%    \bottomrule
% \end{tabular}

\setcounter{section}{0}
\section{Theorem and Corollary }
\textbf{Theorem 1.1 Singular Value Decomposition:} \textit{If $A$ is a real $m$-by-$n$ matrix with rank $p$, then there exist orthogonal matrices:}
$$ U = [u_1|...|u_m]  \in \mathbb{R}^{m*m} \enspace and \enspace V = [v_1|...|v_n]  \in \mathbb{R}^{n*n} $$
such that $$U^TAV = \Sigma = diag(\sigma_1,...\sigma_t)\in \mathbb{R}^{m*n}, t  = min(m,n)$$
retaining only the non-zero portion of the $\Sigma$ matrix, we have
$$A = \sum_{i=1}^p{\sigma_{i}u_iv_i^T} \in \mathbb{R}^{m*n} $$

\textbf{Corollary 1.2 the norm of SVD:} If $A \in \mathbb{R}^{m*n}$ then 
$$||A||_2 = \sigma_1, \   ||A||_F = \sqrt{\sigma_1^2 + \dots + \sigma_p^2},$$

where $p  =min\{m,n\}$.

\section{Preliminaries of norm and eigenvalue}
\textbf{Linear algebra Preliminaries}

\textbf{Lemma 2.1 Frobenius norm of matrix products:}  \textit{For ant matrix $A$ and $B$ we have}

\begin{equation}
\label{lambdaleq}
\begin{split}
||AB||_F\leq||A||_2||B||_F
\end{split}
\end{equation}

\textit{Proof.} Note that:
\begin{equation}
\label{lambdaleq}
\begin{split}
||AB||_F^2 &= \sum_j||(AB)_j||_2^2\leq\sum_j||AB_j||_2^2 \\
&\leq||A||_2^2\sum_j||B_j||_2^2 = ||A||_2||B||_F
\end{split}
\end{equation}

\textbf{Lemma 2.2 Matrix eigenvalue inequality:} \textit{If $A \in \mathbb R^{m*n}$ and $E \in \mathbb R^{m*n}$, then}
\begin{equation}
\label{lambdaleq}
\begin{split}
\sigma_{max}(A+E) \leq \sigma_{max}(A) + ||E||_2 \\
\sigma_{min}(A+E) \leq \sigma_{min}(A) - ||E||_2
\end{split}
\end{equation}

\textit{Proof.} Using Corollary 1.2 it is easy to show that
\begin{equation}
\label{lambdaleq}
\begin{split}
\sigma_{min}(A) \cdot ||x||_2 \leq ||Ax||_2 \leq \sigma_{max}(A) \cdot||x||_2.
\end{split}
\end{equation}
The required inequalities follow from this result.

If a column is added to a matrix, then the largest singular value increases and the smallest singular value decreases.

\textbf{Lemma 2.4 Minkowski inequality:} \textit{If $A \in \mathbb R^{m*n}, B\in \mathbb R^{m*n}$. We have}
\begin{equation}
\label{Minkowski}
\begin{split}
||A+B||\leq||A||+||B||
\end{split}
\end{equation}

\textbf{Lemma 2.5 Cauchy-Schwarz inequality:} \textit{If $A \in \mathbb R^{m*k}, B\in \mathbb R^{k*n}$. We have}
\begin{equation}
\label{Cauchy}
\begin{split}
||AB||\leq \Vert A\Vert \Vert B \Vert
\end{split}
\end{equation}

\textbf{Lemma 2.6 Frobenius norm inequality:} \textit{If $A \in \mathbb R^{m*n}$ and $ (m \ge n)$, has SVD decomposition $A = \sum_{i=1}^n{\sigma_{i}u_iv_i^T}$,We have}
\begin{equation}
\label{Frobenius_inequality}
\begin{split}
 A  \leq \sqrt{n}\sigma_1
\end{split}
\end{equation}

\section{Approximation SVD error analysis }
\textbf{Theorem 3.1 The matrix low-rank approximation using SVD (Singular Value Decomposition):}  \textit{If $A \in \mathbb R^{m*n}$ with rank $p$ has SVD decomposition, take the first \( r(r<p) \) terms to approximate matrix \( A \), then \( A \) can be represented as}

\[ A_r \approx \sum_{i=1}^{r} \sigma_i u_i v_i^T \]

\textit{Then, the factorization satisfies, where $\mathbb E$ denotes expectation with respect to the random test matrix.:}

\begin{equation}
\label{svdfinq}
\begin{split}
\mathbb E||A - A_r||_F \leq \sigma_{k+1}\sqrt{p-r} 
\end{split}
\end{equation}
\textit{Proof.} Due to 
\begin{equation}
\label{lambdaleq}
\begin{split}
||A - A_r||_F &= ||\sum_{i=1}^{p} \sigma_i u_i v_i^T -  \sum_{i=1}^{r} \sigma_i u_i v_i^T||_F \\    
            &= ||\sum_{i=r+1}^{p} \sigma_i u_i v_i^T||_F = \sqrt{\sigma_{r+1}^2 + \dots + \sigma_{p+1}^2}\\
            &\leq  \sigma_{r+1}\sqrt{p-r}
\end{split}
\end{equation}

\textbf{Theorem 3.2 The matrix low-rank approximation using RSVD (Randomized Singular Value Decomposition):}\textit{ If $A \in \mathbb R^{m*n}$. Select an exponent $q$(A parameter in Rsvd) and a target $k$ of singular vectors, where $2\leq k \leq 0.5*min\{m,n\}$. Excute the Randomized SVD algorithm to obtain a rank-$2k$ factorization $A_r = U\Sigma V^{*}. $ Then}
\begin{equation}
\label{rsvderror}
\begin{split}
\mathbb E \Vert A-A_r\Vert \leq \left[ 1+4\sqrt{\frac{2min\{m,n\}}{r-1} }\right ]^{1/(2q+1)}\sigma_{r+1}.
\end{split}
\end{equation}
Where $\mathbb E$ denotes expectation with respect to the random test matrix. The proof of this theorem is complex; for specifics, one may refer to \cite{halko2011finding}.

\section{Analysis of Errors in Quantized operation}
In this section, we use Quantization in $B$ bit and $Quant(*),\widetilde{Quant(*)}$ respectively represent quantization and dequantization.

\textbf{Theorem 4.1 Quantization error of uniformly dithered scalar quantizer:}
\textit{For a scalar $x \in \left [-M, +M \right ]$, denote the quantization error of direct scalar quantizer with a bit-budget of $d$ bits as $\epsilon = Quant(x) - x$. 
The following result further characterizes its mean and variance of this error.}

\begin{equation}
\label{direct_quant_E}
\begin{split}
\mathbb E(\epsilon) = 0\ and\ Var(\epsilon) \leq \frac{M^2}{(2^{d-1}-1)^2}
\end{split}
\end{equation}

\textbf{Theorem 4.2 Quantization error analysis of matrix:} \textit{If $A \in \mathbb R^{m*n}$, and we have $B$ bits quantized martrix $\widetilde{A} = \widetilde{Quant(Quant(A))}'$, denote $\frac{2^{d-1}-1}{M}$ as $\lambda$. Then we have}
\begin{equation}
\label{quantE}
\begin{split}
\Vert A-\widetilde{A} \Vert_F \leq \sqrt{mn}\lambda^{-1} \ and \ \mathbb E(\Vert A-\widetilde{A}\Vert_F) = 0
\end{split}
\end{equation}
\textit{Proof.} 
\textbf{(1).} For the first tern, due to the quantization operation, the error introduced for each element \( a_{ij} \) in matrix \( A \) can be represented as:
\begin{eqnarray}    \label{countAint}
    ra^{Int}_{ij}= \lfloor\lambda *a^{Fp}_{ij}\rfloor - \lambda *a^{Fp}_{ij} 
\end{eqnarray}
where $\lfloor.\rfloor$is the integer type casting(We use downward rounding in our formulas to make it easier to understand). Then
\begin{eqnarray}    \label{eq5}
    ra^{Fp}_{ij}&=&TypeCast((\lfloor\lambda *a^{Fp}_{ij}\rfloor-\lambda *a^{Fp}_{ij})
      *\lambda^{-1},Float) \nonumber \\
      ~&=& TypeCast((ra^{Int}_{ij}*\lambda^{-1},Float)
\end{eqnarray}

Since $ra^{Int}_{ij}\leq1$(Due to the rounding operation, the error in integer will not exceed 1).
Thus 
\begin{equation}
\label{lambdaleq}
\begin{split}
ra^{Fp}_{ij}=TypeCast((ra^{Int}_{ij}*\lambda^{-1},Float)\leq \lambda^{-1}
\end{split}
\end{equation}
The errors \( ra_{ij} \) for each element in matrix \( A \) can be concatenated to form an matrix \( R_A \in \mathbb R^{m*n}\). Then we have 
\begin{equation}
\label{lambdaleq}
\begin{split}
\Vert A-\widetilde{A} \Vert_F &= \Vert R_A \Vert_F= \sqrt{\sum_{i=1}^m\sum_{j=1}^n {ra_{ij}^{Fp}}^2}\\
        &\leq \sqrt{mn}\lambda^{-1}
\end{split}
\end{equation}
\textbf{(2).}For the second term, due to the quantization operation is performed on the entire matrix, we can expand the matrix $A$ into vector. And then we can prove it by applying the previous Theorem.\ref{direct_quant_E} to easily proof.

\textbf{Theorem 4.3 error analysis of Quantized matrix multiplication:} 
\textit{If $C=A*B, \in \mathbb R^{m*n}$ where $A \in \mathbb R^{m*k}, B\in \mathbb R^{k*n}$. Assuming that matrices \( A \) and \( B \) follow the same distribution and $m\leq n \leq k$. Denote the eigenvalue of matrix $A,B$ is $\sigma,\gamma$. For quantized matrix multiplication, $\widetilde{C}=\widetilde{A}*\widetilde{B}$ with $\lambda_1,\lambda_2$. We have}
\begin{equation}
\label{er_quant}
\begin{split}
\mathbb E(\Vert C-\widetilde{C} \Vert_F) &\leq k(\sigma_1\lambda_2^{-1}\sqrt{n}+\gamma_1\lambda_1^{-1}\sqrt{m} \\
    &+\lambda_1^{-1}\lambda_2^{-1}\sqrt{mn})
\end{split}
\end{equation}
\textit{Proof.}
\begin{equation}
\label{lambdaleq}
\begin{split}
\mathbb E(\Vert C-\widetilde{C} \Vert_F)  &= E(\Vert AB-\widetilde{A}\widetilde{B} \Vert_F) \\
        &= E(\Vert (\widetilde{A}+R_A)(\widetilde{B}+R_B)-\widetilde{A}\widetilde{B} \Vert_F) \\
        &= E(\Vert R_A\widetilde{B}+\widetilde{A}R_B+ R_A R_B \Vert_F) 
\end{split}
\end{equation}
Subsequently, by the Minkowski and Cauchy-Schwarz inequality presented in lemma \ref{Cauchy},\ref{Minkowski}
\begin{equation}
\label{lambdaleq}
\begin{split}
\mathbb E &= \mathbb E(\Vert R_A\widetilde{B}+\widetilde{A}R_B+ R_A R_B \Vert_F) \\
    &\leq  \mathbb E(\Vert R_A\widetilde{B} \Vert_F+\Vert \widetilde{A}R_B \Vert_F+ \Vert R_A R_B \Vert_F) \\
    &\leq \mathbb E(\Vert R_A\Vert_F \Vert\widetilde{B}\Vert_F+\Vert \widetilde{A}\Vert_F \Vert R_B \Vert_F + \Vert R_A \Vert_F \Vert R_B \Vert_F)
\end{split}
\end{equation}

Due to previous Theorem.\ref{quantE}. The $\Vert R_A \Vert_F= \leq \mathbb E(\Vert A-\widetilde{A} \Vert_F) \leq \sqrt{mk}\lambda_1^{-1}$, the same way $\Vert R_B \Vert_F \leq \sqrt{kn}\lambda_2^{-1}$. And due to the quantization operation, the maximum absolute value of the quantized matrix will not exceed the maximum absolute value of the original matrix which is $R$, from lemma.(\ref{Frobenius_inequality}) then have 
\begin{equation}
\label{lambdaleq}
\begin{split}
\Vert \widetilde{A} \Vert_F \leq \Vert A \Vert_F = \sqrt{\sum_{i=1}^k {\sigma_{i}}^2} \leq \sigma_1 \sqrt{k}\\
\Vert \widetilde{B} \Vert_F \leq \Vert B \Vert_F = \sqrt{\sum_{i=1}^k {\gamma_{i}}^2} \leq \gamma_1\sqrt{k}\\
\end{split}
\end{equation}
Then incorporating these components into the aforementioned inequality
\begin{equation}
\label{lambdaleq}
\begin{split}
\mathbb E(\Vert C-\widetilde{C} \Vert_F)  &\leq \sigma_1k\lambda_2^{-1}\sqrt{n} +  \gamma_1k\lambda_+1^{-1}\sqrt{m} \\
    &+ k\lambda_1^{-1}\lambda_2^{-1}\sqrt{mn}\\
    &= k(\sigma_1\lambda_2^{-1}\sqrt{n}+\gamma_1\lambda_1^{-1}\sqrt{m}+\lambda_1^{-1}\lambda_2^{-1}\sqrt{mn})
\end{split}
\end{equation}

\begin{algorithm}
    \label{quantization_svd}
  \SetAlgoLined
\KwData{$A$(Intput matrix);$k$(ranf of low rank SVD matrix);$d_1,d_2'$(bit-budget for $\widetilde{U_k}$ and $V_k$)}
  
  \KwResult{$\widetilde{U_k}',V_k'$(The matrix after quantization)}

    \tcc{Computes SVD approximate matrix of $A$}
    $\{U_{k},\Sigma_{k},V_{k}^{T}) \}\leftarrow SVD(A,k)$\;

    \tcc{Calculate the product of the diagonal matrix $\widetilde{U_k}$}
    $\widetilde{U_k}\leftarrow U_{k} \cdot \Sigma_{k}$\;    

    \tcc{Calculate the quantization svd approximating matrix}
    $\{\widetilde{U_k}',V_k'\}  \leftarrow \widetilde{Quant}(Quant(\{\widetilde{U_k},V_k\},\{d_1,d_2\}),T_F)$\;

    \Return  $\widetilde{U_k}', V_k'$\;
    
  \caption{Algorithms of quantization svd approximating matrix A.}
\end{algorithm}

\textbf{Theorem 4.4 error analysis of quantization SVD\cite{saha2023matrix}:} 
\textit{For a matrix $A \in \mathbb R^{m*n}$ and $m\leq n$ has SVD approximation matrix $A' = U\Sigma V^T$, followed by quantizing the singular vectors corresponding to the top-k singular values, scaled by the singular values. Let $\widetilde U = U\Sigma$ and the first k column as $\widetilde{U_k}$.  Assuming that \( d_1 \) and \( d_2 \) are used to quantize matrices \( \widetilde{U_k} \) and \( V_k \) respectively, we obtain the quantized matrices \( \widetilde{U_k}', V_k'\) with \( \lambda_1 \) and \( \lambda_2 \). Then we have the quantization svd approximating matrix $\widetilde {A'}$, the pseudocode is provided in Algorithm.\ref{quantization_svd}.}

\textit{Then we have}
\begin{equation}
\label{E_qsvd_proof0}
\begin{split}
\mathbb E(\Vert \widetilde{U_k}'V_k'^T - A \Vert_F^2) &\leq  \sigma_{k+1}^2(n-r) + mk\frac{\sigma_1^2}{(2^{d_1-1}-1)^2} \\
    &+ nk\sigma_1^2\frac{1}{(2^{d_2-1}-1)^2} \\
    &+ mnk\frac{\sigma_1^2}{(2^{d_1-1}-1)^2}\frac{1}{(2^{d_2-1}-1)^2}
\end{split}
\end{equation}

\textit{Proof.}
Denote $R_U = \widetilde{U_k}' - \widetilde{U_k}$ and $R_V = V_k'-V_k$. The approximation error is given by
\begin{equation}
\label{E_qsvd_proof1}
\begin{split}
\mathbb E&(\Vert \widetilde{U_k}'V_k'^T - A \Vert_F^2) \\
        &= \mathbb E(\Vert \widetilde{U_k}'(V_k^T+R_V^T) - A \Vert_F^2) \\
        &\stackrel{(i)}{=}  \underbrace{\mathbb E(\Vert \widetilde{U_k}'V_k^T - A \Vert_F^2)}_{T_1}  + \underbrace{\mathbb E(\Vert \widetilde{U_k}'R_V^T\Vert_F^2)}_{T_2}
\end{split}
\end{equation}

where the cross term disappears in $(i)$ because by Theorem.\ref{quantE} we can get $\mathbb E(R_V)=0$ . Then
$$\mathbb E(Tr((\widetilde{U_k}'V_k^T)))^T\widetilde{U_k}\mathbb E(R_V)) = 0$$.

\textbf{In $T_1$:} The first term in (\ref{E_qsvd_proof1}) can be upper bounded as
\begin{equation}
\label{E_qsvd_proof2}
\begin{split}
\mathbb E(\Vert \widetilde{U_k}'V_k^T - A \Vert_F^2) &= \mathbb E(\Vert (\widetilde{U_k}+R_U)V_k^T - A \Vert_F^2) \\
        &\stackrel{(ii)}{=} \Vert A_k - A \Vert_F^2 + \mathbb E(\Vert R_UV_k^T \Vert_F^2)
\end{split}
\end{equation}
where the cross term disappears in $(ii)$ because also $\mathbb E(R_U)=0$. Noting that each column of matrix \( V_k \) is an orthogonal vector then $V_kV_k^T=I$. The second trem in (\ref{E_qsvd_proof2}) is 
\begin{equation}
\label{E_qsvd_proof3}
\begin{split}
\mathbb E(\Vert R_UV_k^T \Vert_F^2) &= \mathbb E( Trace(R_U V_k^T V_k R_U^T)) \\
        &= Trace(\mathbb E(R_U R_U^T)) \stackrel{(iii)}{\leq} mk\lambda_1^{-2}
\end{split}
\end{equation}
Here (iii) follows the Theorem.\ref{quantE}, then for matrix $R_U$
$$\mathbb E(Trace(R_UR_U^T)) = \mathbb E( \sum_{i=1}^m\sum_{j=1}^k r_{Uij}^2 ) \leq mnVar(r_U) = mk\lambda_1^{-2}$$

Then (\ref{E_qsvd_proof2}) can be upper bounded as
\begin{equation}
\label{E_qsvd_proof4}
\begin{split}
\mathbb E(\Vert \widetilde{U_k}'V_k^T - A \Vert_F^2) \leq \Vert A_k - A \Vert_F^2 + mk\lambda_1^{-2}
\end{split}
\end{equation}

\textbf{In $T_2$:} The second term in (\ref{E_qsvd_proof1}) can be upper bounded as
\begin{equation}
\label{E_qsvd_proof5}
\begin{split}
\mathbb E(\Vert \widetilde{U_k}'R_V^T\Vert_F^2) &= \mathbb E(\Vert (\widetilde{U_k}+R_U)R_V^T\Vert_F^2)\\ &= \mathbb E\Vert \widetilde{U_k}R_V^T \Vert_F^2 + E\Vert R_UR_V^T \Vert_F^2
\end{split}
\end{equation}

where the cross term disappears due to $\mathbb E(R_U)=0$ and can proof as the term in $T_1$. Then the first term in (\ref{E_qsvd_proof5}) is
\begin{equation}
\label{E_qsvd_proof5}
\begin{split}
\mathbb E\Vert \widetilde{U_k}R_V^T \Vert_F^2 &= \mathbb E(Trace(R_U^T\widetilde{U_k}^T\widetilde{U_k}R_U)) \\
    &= \mathbb E(Trace(\widetilde{U_k}^T\widetilde{U_k}R_UR_U^T)) \\
    &= \mathbb E(Trace(\widetilde{U_k}^T\widetilde{U_k}))\mathbb E(R_UR_U^T)_{i=j}\\
    &\stackrel{(iv)}{\leq} n\lambda_2^{-2}\mathbb E(Trace(\widetilde{U_k}^T\widetilde{U_k}))\\
    &\stackrel{(v)}{\leq} n\lambda_2^{-2}\sum_{i=1}^k{\sigma_i^2}\leq nk\sigma_1^2\lambda_2^{-2}
\end{split}
\end{equation}
Here $(iv)$ follows the Theorem.\ref{quantE}, that
$$
    \mathbb E(R_UR_U^T)_{i=j} = avg(\sum_{i=0}^n r_{Uij}^2) = n\lambda_2^{-2}
$$
And $(v)$ follows that the matrix $\widetilde{U_k} = U\Sigma = \sum_{i=0}^k\sigma_i u_i$. Then 
\begin{equation}
\label{E_qsvd_proof6}
\begin{split}
\widetilde{U_k}^T\widetilde{U_k}) &= \sum_{i=0}^k(\sigma_i u_i)(\sigma_i u_i)^T = \sum_{i=0}^k \sigma_i^2 \Vert u_i \Vert_2^2 \\
    &=  \sum_{i=0}^k \sigma_i^2 \leq k\sigma_1^2
\end{split}
\end{equation}

The second term in \ref{E_qsvd_proof5} is 
\begin{equation}
\label{E_qsvd_proof7}
\begin{split}
    \mathbb E\Vert R_UR_V^T \Vert_F^2 &= \mathbb E\Vert Trace(R_U^TR_VTR_V^TR_U) \Vert_F^2 \\
        &= \mathbb E\Vert  Trace(R_VR_V^T) \mathbb E(R_UR_U^T)_{i=j} \Vert_F^2\\
        &\stackrel{(vi)}{\leq} mnk\lambda_1^{-2}\lambda_2^{2}
\end{split}
\end{equation}
The proof to \textbf{(vi)} is same to \textbf{(iv)} and \textbf{(iii)}.

Summarizing $T_2$ and $T_3$ we have 
\begin{equation}
\label{E_qsvd_proof8}
\begin{split}
\mathbb E(\Vert \widetilde{U_k}'V_k'^T - A \Vert_F^2) &\leq  \Vert A_k-A \Vert_F^2 + mk\lambda_1^{-2} \\
    &+ nk\sigma_1^2\lambda_2^{-2} + mnk\lambda_1^{-2}\lambda_2^{-2}
\end{split}
\end{equation}

From Theorem.\ref{svdfinq} $$\mathbb E(\Vert A - A_k\Vert_F)\leq \sigma_{k+1}\sqrt{n-r}$$

Noting that matrix $\widetilde{U_k} = U_k*\Sigma_k$ and $V_k,U_k$ are orthogonal matrix. Then the $|max(\widetilde{U_k}_{ij})| \leq \sigma_1$ and $|max(V_{kij})| \leq 1$. Then 
\begin{equation}
\label{lambdaleq}
\begin{split}
    \lambda_1 \ge \frac{2^{d_1-1}-1}{\sigma_1} \ , \ \lambda_2 \ge \frac{2^{d_1-1}-1}{1}
\end{split}
\end{equation}

Incorporate these two terms into (\ref{E_qsvd_proof8}), we have 
\begin{equation}
\label{E_qsvd_proof9}
\begin{split}
\mathbb E(\Vert \widetilde{U_k}'V_k'^T - A \Vert_F^2) &\leq  \sigma_{k+1}^2(n-r) + mk\frac{\sigma_1^2}{(2^{d_1-1}-1)^2} \\
    &+ nk\sigma_1^2\frac{1}{(2^{d_2-1}-1)^2} \\
    &+ mnk\frac{\sigma_1^2}{(2^{d_1-1}-1)^2}\frac{1}{(2^{d_2-1}-1)^2}
\end{split}
\end{equation}

which is (\ref{E_qsvd_proof0}).

\section{Analysis of Errors in Quantized AMM via SVD}
\textit{In the preceding sections, we analyzed the errors associated with quantization operations and the use of SVD approximation on individual matrices. In this section, we will combine the quantized matrix multiplication with the SVD approximation to analyze the ultimate error of the algorithm.
}
\textit{In this section, we focus on quantized AMM $C' = A'B'$. The  the pseudocode is provided in Algorithm\ref{qrsvdgemm}. }

\textit{Uisng rank $r$ as the approximation matrix rank. Similar to the previous section matrix $A \in \mathbb R^{m*k}$ has SVD approximation matrix $A_r = U_r\Sigma_r V_r^T = \widetilde{U_r}V_r$ with quantization bit $d_1,d_2$, matrix $B \in \mathbb R^{k*n}$ has SVD approximation matrix $B_r =  W_{r}\Gamma_{r}Z_{r}^{T} = W_k\widetilde{Z_k}$ with quantization bit $d_4,d_3$. Assume $m\leq n \leq k$}.

\textit{Due to matrix multiplication, the quantization bit $d_2$ must equal to $d_2$. Then we have quantization matrix $\widetilde{U_r},V_r,W_r,\widetilde{Z_r}$ with $d_1,d_2,d_2,d_3$.
So
$${C'_r} = {A'_r}{B'_r} = \widetilde{U_r'}V_r'W_r'\widetilde{Z_r'}$$
\textit{  Denote $2^{d_i-1}-1$} as $D_i$. Then we have}
\begin{equation}
\label{E_qamm_proof0}
\begin{split}
\mathbb E(\Vert {C'_r} - C \Vert_F) \leq \sqrt{\mathbb L_1* \mathbb L_4} + \sqrt{\mathbb L_1* \mathbb L_2}  + \sqrt{\mathbb L_3* \mathbb L_4}  
\end{split}
\end{equation}
\textit{Where $\mathbb L_1 \sim  \mathbb L_4$ are given in below Proof.}

\textit{Proof. }
Denote $R_A =  {A'_r} - A$, $R_B =  {B'_r} - B$. Then 
\textit{Then we have}
\begin{equation}
\label{E_qamm_proof1}
\begin{split}
\mathbb E(\Vert {C'_r} - C \Vert_F) &= \mathbb E(\Vert \widetilde{A'}\widetilde{B'}- AB \Vert_F) \\
        &= \mathbb E(\Vert (A+R_A)(B+R_B)- AB \Vert_F)\\
        &=  \mathbb E(\Vert R_AB+R_BA + R_AR_B \Vert_F)\\
        &\leq \mathbb E(\Vert R_A \Vert_F\Vert B \Vert_F+\Vert R_B \Vert_F\Vert A \Vert_F \\
        &+\Vert R_A \Vert_F\Vert R_B \Vert_F)
\end{split}
\end{equation}
From Theorem.\ref{E_qsvd_proof0}, we have 
\begin{equation}
\label{E_qamm_proof1}
\begin{split}
\mathbb L_1 &= \mathbb E(\Vert R_A \Vert_F^2) = \mathbb E(\Vert \widetilde{U_r'}V_r - A \Vert_F^2) \\
        &\leq \sigma_{r+1}^2(k-r) + mr\frac{\sigma_1^2}{(2^{d_1-1}-1)^2} \\
    &+ kr\sigma_1^2\frac{1}{(2^{d_2-1}-1)^2} \\
    &+ mkr\frac{\sigma_1^2}{(2^{d_1-1}-1)^2}\frac{1}{(2^{d_2-1}-1)^2}
\end{split}
\end{equation}\
Similar derivation applied to \( R_B \) (here $\gamma_1$ is the first eigenvalue of matrix $B$ different from rank $r$)
\begin{equation}
\label{E_qamm_proof2}
\begin{split}
\mathbb L_2 &=\mathbb E(\Vert R_B \Vert_F^2) = \mathbb E(\Vert W_r\widetilde{Z_r'} - B \Vert_F^2) \\
        &\leq \gamma_{r+1}^2(k-r) + kr\frac{\gamma_1^2}{(2^{d_3-1}-1)^2} \\
    &+ nr\gamma_1^2\frac{1}{(2^{d_2-1}-1)^2} \\
    &+ knr\frac{\gamma_1^2}{(2^{d_3-1}-1)^2}\frac{1}{(2^{d_2-1}-1)^2}
\end{split}
\end{equation}

And from rom lemma.(\ref{Frobenius_inequality}) 
\begin{equation}
\label{lambdaleq}
\begin{split}
\mathbb L_3 = \Vert A \Vert_F^2 \leq \sigma_1^2 {k} \ \ , \ \ 
\mathbb L_4 = \Vert B \Vert_F^2 \leq \gamma_1^2{k}
\end{split}
\end{equation}

\section{Analysis of Errors in Quantized AMM via SVD under specific distribution}
\textit{In the previous section, we analyzed the general error case of quantized AMM. However, the form presented contained too many parameters, which was not conducive to intuitive analysis. In this section, we will impose certain restrictions on matrices \( A \) and \( B \), as well as the quantization conditions, to provide a more intuitive form of error.}

Assuming the matrix dimensions are \( m = n = k \), and that matrices \( A \) and \( B \) are identically distributed with a maximum value of \( R \), and let $ \sigma_1 = \gamma_1 $, $\sigma_r = \gamma_r = \sigma_1/\sqrt{r}$. 

Denote $2^{d_i-1}-1$ as $D_i$, $\frac{\sigma_{r+1}^2}{\sigma_1^2}(k-r)$ as $f(r)$ which means under a fixed distribution for matrices, the magnitude of this term depends solely on the value chosen for \( r \).
Then
\begin{equation}
\label{spec_amm_1}
\begin{split}
\mathbb E(\Vert {C'_r} - C \Vert_F) &\leq k\sigma_1^2\sqrt{r} \bigg(\sqrt{f(r)+\frac{1}{D_1}+\frac{1}{D_2}+\frac{K}{D_1D_2}}  \\
        &+ \sqrt{f(r)+\frac{1}{D_2}+\frac{1}{D_3}+\frac{K}{D_2D_3}} \bigg ) \\
        &+kr\sigma_1^2\bigg( \sqrt{f(r)+\frac{1}{D_1}+\frac{1}{D_2}+\frac{K}{D_1D_2}} \bigg)\\
        &\cdot \bigg ( \sqrt{f(r)+\frac{1}{D_2}+\frac{1}{D_3}+\frac{K}{D_2D_3}} \bigg)
\end{split}
\end{equation}
\\\
\textit{Proof.}
Then the first term in (\ref{E_qamm_proof2}) can be upper bounded as.

\begin{equation}
\label{spec_amm_2}
\begin{split}
\mathbb E(\Vert R_A \Vert_F\Vert B \Vert_F) &\leq \Bigg (\gamma_1^2\sigma_{r+1}^2k(k-r) + mrk\frac{\gamma_1^2\sigma_1^2}{D_1} \\
    &+ mrk^2\frac{\gamma_1^2\sigma_1^2}{D_1D_2}\Bigg) ^\frac{1}{2} \\
    &=  \sigma_1\sqrt{kr} \sqrt{f(r)+\frac{1}{D_1}+\frac{1}{D_2}+\frac{K}{D_1D_2}}
\end{split}
\end{equation}

Then the second term in (\ref{E_qamm_proof2}) can be upper bounded as
\begin{equation}
\label{spec_amm_3}
\begin{split}
\mathbb E(\Vert R_B \Vert_F\Vert A \Vert_F) &\leq \Bigg (\gamma_{r+1}^2\sigma_{1}^2k(k-r) + nrk\frac{\gamma_1^2\sigma_1^2}{D_3} \\
    &+ nrk^2\frac{\gamma_1^2\sigma_1^2}{D_2D_3}\Bigg) ^\frac{1}{2} \\
    &=  \sigma_1\sqrt{kr} \sqrt{f(r)+\frac{1}{D_2}+\frac{1}{D_3}+\frac{K}{D_2D_3}}
\end{split}
\end{equation}
By combining the two terms (\ref{E_qamm_proof1})(\ref{E_qamm_proof2}) into (\ref{E_qamm_proof0}), the proof can be established.

\end{document}